\documentclass[11pt,a4paper,oneside,reqno]{amsart}

\usepackage{setspace}

\usepackage[T1]{fontenc}
\usepackage[english]{babel}
\usepackage{lmodern}
\usepackage{microtype}

\usepackage[autostyle]{csquotes}
\MakeAutoQuote{‘}{’}
\MakeOuterQuote{"}

\usepackage{amsmath,amssymb,amscd,amsthm}
\usepackage{mathtools}
\usepackage{stmaryrd}
\usepackage{latexsym}
\usepackage{bbm,amsbsy}
\usepackage{mathrsfs}
\usepackage{nicefrac}
\usepackage{xfrac}
\usepackage{faktor}

\usepackage{graphicx}
\usepackage{epsfig}
\usepackage[all,cmtip]{xy}
\SelectTips{eu}{12}
\usepackage{tikz}
\usepackage{tikz-cd}
\usepackage{caption}
\usepackage{subfigure}

\usepackage{colortbl}
\usepackage{color}
\usepackage{enumitem}
\usepackage{scalerel}
\usepackage{cancel}
\usepackage{systeme}
\usepackage{nomencl}
\usepackage[normalem]{ulem}

\usepackage[colorlinks=true,pagebackref=true]{hyperref}

\makenomenclature

\makeatletter
\newtheorem*{rep@theorem}{\rep@title}
\newcommand{\newreptheorem}[2]{%
\newenvironment{rep#1}[1]{%
 \def\rep@title{#2 \ref{##1}}%
 \begin{rep@theorem}}%
 {\end{rep@theorem}}}
\makeatother

\makeatletter
\newtheorem*{rep@cor}{\rep@title}
\newcommand{\newrepcor}[2]{%
\newenvironment{rep#1}[1]{%
 \def\rep@title{#2 \ref{##1}}%
 \begin{rep@cor}}%
 {\end{rep@cor}}}
\makeatother

\makeatletter
\newtheorem*{rep@prop}{\rep@title}
\newcommand{\newrepprop}[2]{%
\newenvironment{rep#1}[1]{%
 \def\rep@title{#2 \ref{##1}}%
 \begin{rep@prop}}%
 {\end{rep@prop}}}
\makeatother

\newtheorem{cor}{Corollary}[section]   

\newtheorem{theorem}[cor]{Theorem}
\newtheorem{prop}[cor]{Proposition}

\newtheorem{fact}[cor]{Fact}

\newtheorem{remark}[cor]{Remark}

\newtheorem{lemma}[cor]{Lemma}

\newtheorem{corx}{Corollary}

\newtheorem{thmx}[corx]{Theorem}
\newtheorem{propx}[corx]{Proposition}
 
\newrepcor{cor}{Corollary}
\newreptheorem{theorem}{Theorem}
\newrepprop{prop}{Proposition}

\newlist{steps}{enumerate}{1}
\setlist[steps, 1]{itemsep=8pt,leftmargin=0cm,itemindent=.5cm,labelwidth=\itemindent,labelsep=0cm,align=left,label = \textbf{\emph{Step \arabic*}:\,}}

\makeatletter
\newcommand{\myitem}[1]{%
\item[#1]\protected@edef\@currentlabel{#1}%
}
\makeatother

\newcommand{\R}{\mathbb{R}}

\newcommand{\SL}{\operatorname{SL}}

\newcommand{\SO}{\operatorname{SO}}

\newcommand{\C}{\mathbb{C}}

\def\C{\mathbb{C}}
\def\R{\mathbb{R}}

\def\N{\mathbb{N}}

\def\Aff{\operatorname{Aff}}

\def\GL{{\operatorname{GL}}}

\def\SO{{\operatorname{SO}}}
 
\def\SU{{\operatorname{SU}}} 
\def\SL{{\operatorname{SL}}}

\def\Stab{{\operatorname{Stab}} }

\def\p{{\mathfrak{p}}}

\def\k{{\mathfrak{k}}}

\def\a{{\mathfrak{a}}}
\def\g{{\mathfrak{g}}}

\def\n{{\mathfrak{n}}}

\title[Kleinian affine manifolds]{On divisible domains of the flat affine space with parallel volume}
\author[Farid Diaf, Blandine Galiay, Malek Hanounah]{Farid Diaf \and Blandine Galiay \and Malek Hanounah}

\address{Farid Diaf: Institut de Recherche Mathématique Avancée, UMR 7501, Université de Strasbourg et CNRS, 7 rue René Descartes, 67000 Strasbourg, France.} 
\email{f.diaf@unistra.fr}

\address{Blandine Galiay: CNRS et Institut \'Elie Cartan de Lorraine, Faculté des Sciences et Technologies, Campus, Boulevard des Aiguillettes, 54506 Vandœuvre-lès-Nancy, France.}
\email{blandine.galiay@univ-lorraine.fr}

\address{Malek Hanounah: MPI für Mathematik in den Naturwissenschaften,
Inselstraße 22, 04103 Leipzig
Germany} 
\email{malek.hanounah@mis.mpg.de}

\begin{document}\raggedbottom

\setcounter{secnumdepth}{3}
\setcounter{tocdepth}{2}

\begin{abstract}
We consider a domain $\Omega\subset\R^d$ divisible by a
discrete subgroup $\Gamma$ of the unimodular affine group
$\SL_d(\R)\ltimes\R^d$.
We prove that $\Omega=\R^d$ if $\Gamma$ is virtually solvable
or if its linear holonomy is virtually contained in a
semisimple Lie subgroup of real rank one.
In particular, we establish the Kleinian Markus conjecture
in both cases.
\end{abstract}

\maketitle

\tableofcontents

\section{Introduction}

The study of geometric structures on manifolds is an important topic in geometry and topology. Among these, \emph{affine structures} have received considerable attention and remain an active area of research; see, for instance, \cite{milnor1977fundamental,goldman1986affine,carriere1989autour,
DGK_margulis}.

An affine manifold is a manifold endowed with an
$(\Aff(\R^d),\R^d)$-structure, where
$\Aff(\R^d)=\GL_d(\R)\ltimes\R^d$ is the group of
affine transformations of $\R^d$.
One of the main open questions in affine geometry is
the \emph{Markus conjecture}, which asserts that every
closed affine manifold with parallel volume is complete.
Here, having parallel volume means that the linear
holonomy is contained in $\SL_d(\R)$.

In this paper, we study the Markus conjecture in the
\emph{Kleinian setting}. By a Kleinian affine manifold,
we mean a quotient $M=\Gamma\backslash\Omega$, where
$\Omega\subset\R^d$ is a domain, that is, a nonempty
connected open subset, and
$\Gamma\leq\Aff(\R^d)$ is a discrete subgroup acting
freely, properly discontinuously, and cocompactly on
$\Omega$. Such a domain is called \emph{divisible}.
A Kleinian affine manifold is complete precisely when
$\Omega=\R^d$.

The \emph{Kleinian Markus conjecture} asserts that
if a domain $\Omega\subset\R^d$ is divisible by a
discrete subgroup
$\Gamma\leq\SL_d(\R)\ltimes\R^d$, then
$\Omega=\R^d$. The conjecture remains open in general, although
several partial results are known.

\subsection{State of the art}

The Markus conjecture has been established under various
algebraic and geometric restrictions.

On the algebraic side, completeness is known when the
affine holonomy is virtually abelian
\cite{smillie1977affinely} or virtually nilpotent
\cite{fried1981affine}.
Goldman and Hirsch \cite{goldman1986affine} also
established completeness for solvable affine holonomy
under an additional bound on its solvable rank.

On the geometric side, an important notion introduced
by Carri\`ere is that of \emph{discompacity}.
Roughly speaking, the discompacity of a linear group
measures the maximal number of independent directions
along which Euclidean balls can collapse under the
action of sequences of elements of the group.

Carri\`ere \cite{carriere1989autour} proved the Markus
conjecture when the linear holonomy has discompacity
at most one. This includes the case of flat Lorentzian
manifolds, since standard representation of $\SO(n,1)$ on $\mathbb{R}^{n+1}$ has discompacity one. Tholozan \cite{tholozan2015completude} obtained further completeness results in the Kleinian setting for
linear holonomy contained in $\SU(n,1)$ and
$\SO(3,\C)$, acting through their standard complex
representations on $\mathbb{C}^{n+1}$ and $\mathbb{C}^3$ respectively. These groups have real rank one,
and their standard representations have "complex"
discompacity one.

More recently, the Kleinian Markus conjecture was
established by the authors for linear holonomy contained in
$\SO(2,2)$ \cite{diaf2026completeness}, whose standard
representation on $\mathbb{R}^4$ has discompacity two.
Other results concern non-Zariski-dense linear
holonomy in $\SL_3(\R)$ \cite{hanounahmarkus3d}. Results of a different nature have been obtained
by imposing geometric conditions on the domain
rather than on the linear holonomy. In particular,
the conjecture is known for convex divisible domains
in dimensions at most five \cite{Markus_convex}.

\subsection{Main results}

Our first main result establishes the Kleinian Markus
conjecture when the affine holonomy is virtually solvable.

\begin{thmx}\label{thm:solvable-complete}
Let $\Omega\subset\R^d$ be a domain, and let
$\Gamma\leq\SL_d(\R)\ltimes\R^d$ be a discrete,
virtually solvable subgroup acting freely, properly
discontinuously, and cocompactly on $\Omega$.
Then $\Omega=\R^d$.
\end{thmx}

Goldman and Hirsch \cite{goldman1986affine} established
a completeness result for closed affine manifolds with
solvable holonomy under a restriction on the solvable rank.
Our theorem removes this restriction in the Kleinian
setting, allowing virtually solvable holonomy groups
of arbitrary solvable rank.

Solvable groups also play an important role in affine geometry.
For instance, the Auslander conjecture \cite{auslander1964structure}
asserts that the fundamental group of every closed complete
affine manifold is virtually solvable. Our theorem goes in the other direction in the
Kleinian setting: virtual solvability, together with parallel
volume, implies completeness.

Our second main result concerns linear holonomy
contained in a semisimple Lie group of real rank one.

\begin{thmx}\label{thm:rank1-complete}
Let $G$ be a connected real semisimple Lie group
of real rank one, and let
$$
\rho:G\longrightarrow\GL_d(\R)
$$
be a faithful representation.
Let $\Omega\subset\R^d$ be a domain, and let
$\Gamma\leq\rho(G)\ltimes\R^d$ be a discrete
subgroup acting freely, properly discontinuously,
and cocompactly on $\Omega$.
Then $\Omega=\R^d$.
\end{thmx}

Note that since $G$ is connected and semisimple, the image of
$\rho$ is contained in $\SL_d(\R)$. Our theorem extends in the Kleinian setting the results of Carrière \cite{carriere1989autour} and Tholozan \cite{tholozan2015completude} for the standard
representations of $\SO(n,1)$, $\SU(n,1)$, and
$\SO(3,\C)$ to arbitrary faithful representations of
real-rank-one semisimple Lie groups. Recall that Carrière's result also aplies to non-Kleinian affine manifolds. The rank-one case of the Markus conjecture was also raised by Kapovich in \cite[Conjecture 2.12]{kapovichBIRS2019}.

An important feature of Theorem \ref{thm:rank1-complete}
is that the discompacity of the linear holonomy
can be arbitrarily large, even though it is contained in a semisimple Lie group of real rank one. For example, consider the diagonal embedding
$$
\rho_n:\SL_2(\R)\longrightarrow\SL_{2n}(\R),
\qquad
\rho_n(g)=\operatorname{diag}(g,\ldots,g).
$$
The matrices $\rho_n(\operatorname{diag}(e^t,e^{-t}))$
have $n$ contracting directions as $t\to+\infty$. In fact, one can show that the discompacity of $\rho_n(\SL_2(\R))$ is exactly $n$.
Thus, our theorem applies to linear groups of arbitrarily
large discompacity, going beyond the discompacity-one
setting of Carri\`ere's theorem \cite{carriere1989autour} and the discompacity-two setting of the $\SO(2,2)$ case \cite{diaf2026completeness}.

\subsection{Outline of the proof}

We briefly describe the main ideas behind the proofs
of our two theorems.

The key ingredient in the proof of
Theorem \ref{thm:solvable-complete} is the following
general dynamical result.

\begin{propx}\label{prop:dynamic}
Let $G$ be a Lie group acting on a manifold $X$,
and let $\Gamma$ be a lattice in $G$.
Assume that $\Gamma$ acts properly discontinuously
on a domain $\Omega\subset X$.
Then, for every $x\in\Omega$, the identity component
of the stabilizer $G_x$ is compact.
\end{propx}

When $\Gamma$ is virtually solvable, we may pass to
a finite-index solvable subgroup and use the existence
of a \textit{syndetic hull} $S$ containing $\Gamma$ as a uniform lattice.
Applying Proposition \ref{prop:dynamic} with $G=S$
and $X=\R^d$, we obtain that the identity component
of the stabilizer $S_x$ is compact for every $x\in\Omega$.
This gives a bound on the cohomological dimension of
$\Gamma$, and completeness follows from a criterion
of Goldman and Hirsch \cite{goldman1986affine}.

We now turn to the rank-one result. The proof of
Theorem \ref{thm:rank1-complete} relies on the following
general reduction result. To state it, let $G\leq\SL_d(\R)$ be a semisimple Lie subgroup, and fix a Cartan decomposition $G=KA^+K$
and a $K$-invariant Euclidean norm $|\cdot|$ on $\R^d$.
For a subgroup $G'\leq G$, let $A^+(G')$ denote the set
of Cartan projections of elements of $G'$.
We say that a subset $X\subset G$ \emph{does not contract}
a subspace $E\subset\R^d$ with respect to $|\cdot|$
if $E$ is invariant under $X$ and
$|g\cdot v|\geq |v|$ for all $g\in X$ and $v\in E$.

\begin{propx}\label{mainprop}
Let $\Omega\subsetneq\R^d$ be a domain, and let
$\Gamma\leq\SL_d(\R)\ltimes\R^d$ be a discrete
subgroup acting freely, properly discontinuously,
and cocompactly on $\Omega$.

Suppose that the linear part $L(\Gamma)$ is contained
in a semisimple Lie subgroup $G\leq\SL_d(\R)$
satisfying the following conditions:
\begin{enumerate}
    \item There exists an integer $p\geq d/2$ and a
    $G$-invariant orbit
    $\mathcal F\subset\operatorname{Gr}_p(\R^d)$ 
    whose closure $\overline{\mathcal{F}}$ has the following property: any distinct subspaces $E$, $F\in \overline{\mathcal{F}}$ are transverse, that is, $E+F=\mathbb{R}^d$.

    \item There exists $E\in\mathcal F$ such that
    $A^+(L(\Gamma))$ does not contract $E$.
\end{enumerate}

Then there exists an $L(\Gamma)$-invariant subspace
$E'\in\overline{\mathcal F}$ such that $\Omega$
is foliated by affine subspaces parallel to $E'$.
\end{propx}

In several classical examples, the family $\mathcal F$
has a simple geometric description. For the standard
representation of $\SO_0(n,1)$ on $\R^{n+1}$,
it consists of the coisotropic hyperplanes with
respect to the preserved Lorentzian form.
Similarly, in the cases of $\SU(n,1)$ and
$\SO(3,\C)$ considered by Tholozan
\cite{tholozan2015completude}, the family
$\mathcal F$ consists of coisotropic complex
hyperplanes with respect to the corresponding
Hermitian or complex quadratic form respectively.

In Proposition \ref{thm_rank_one}, we show that the hypotheses of
Proposition \ref{mainprop} are satisfied for
arbitrary faithful finite-dimensional representations
of real-rank-one semisimple Lie groups.
The essential ingredient is a transversality
phenomenon for the orbit of the subspace of
nonnegative weights established in Lemma \ref{lemma:transversality rank 1}. This extends the geometric
picture of coisotropic hyperplanes described above
to representations for which the corresponding
subspaces need not be hyperplanes.

As a consequence, we will show that the linear holonomy preserves a
subspace whose stabilizer is a    parabolic
subgroup. Since    parabolic subgroups are
amenable, a result of Milnor \cite{milnor1977fundamental}
implies that $\Gamma$ is virtually polycyclic,
and hence virtually solvable. The conclusion
then follows from Theorem \ref{thm:solvable-complete}.

\subsection{Organization of the paper}

In Section \ref{Section: Dynamical lemma}, we prove
the dynamical result of Proposition \ref{prop:dynamic}
and use it to establish Theorem \ref{thm:solvable-complete}.
Section \ref{sec:rank1} is devoted to the rank-one results.
We first prove the general geometric reduction result
of Proposition \ref{mainprop}, then establish the
transversality property for representations of
real-rank-one semisimple Lie groups, and finally
prove Theorem \ref{thm:rank1-complete}.

\vspace{1mm}

\textbf{Use of AI:} The first author used an LLM while developing
Lemma \ref{lemma:dynamic}. The interaction began with discrete solvable groups whose linear parts are upper triangular matrices, with the aim of understanding the relation between their syndetic hulls and divisible domains. At one point, the LLM suggested using the Poincar\'e recurrence theorem. The author then explored how far this method could be extended through several rounds of precise instructions and wrote the final statement and proof of the lemma. The lemma is related in spirit to ideas appearing in the works of Kobayashi \cite{Kobayashi1989ProperAction} and Frances \cite{frances2005lorentzian}. Aside from this lemma, all mathematical arguments and the overall strategy of the paper were developed by the authors.

\section{The solvable case}\label{Section: Dynamical lemma}

Th goal of this section is to prove Theorem \ref{thm:solvable-complete}.

\subsection{Dynamical results}

The following lemma will be the key to show Proposition \ref{prop:dynamic}.

\begin{lemma}\label{lemma:dynamic}
Let $G$ be a Lie group acting on a manifold $X$. Let $\Gamma$ be a lattice
in $G$, and assume that $\Gamma$ acts properly discontinuously on a domain
$\Omega \subset X$. For $x \in \Omega$, let $G_x$ denote the stabilizer of $x$
in $G$. Then, for every $h \in G_x$, the sequence $(h^n)_{n \geq 1}$ admits
a convergent subsequence.
\end{lemma}
Note that the subtlety in this lemma is that $\Omega$ need not be invariant under $G$. Nevertheless, we will see that the lattice assumption allows us to obtain information about the stabilizer $G_x$ of a point $x\in \Omega$.   

It is a classical fact that, when $\Gamma\leq G$ is a uniform lattice which acts properly discontinuously and cocompactly on the whole space $X$, then the action of $G$ on $X$ is proper
    (see, for instance, \cite[Lemma 2.3]{Kobayashi1989ProperAction}).
    Therefore, the conclusion of the last proposition is immediate in this case, since the stabilizers are compact.
 Another observation, coming from affine geometry, is that if a solvable
    subgroup $\Gamma \leq \mathrm{Aff}(\mathbb{R}^d)$ acts properly discontinuously
    and cocompactly on $\mathbb{R}^d$, then the stabilizers $H_x$ in the
    syndetic hull $H$ of $\Gamma$ are finite \cite[Theorem 1.4]{Crysta_FGH}. This observation is one of the
    motivations for considering the previous proposition, and it will also be
    useful later in the proof of Theorem \ref{thm:solvable-complete}.

Before proving the proposition, recall that since $\Gamma$ is a lattice,
the quotient $\Gamma\backslash G$ carries a finite $G$-invariant measure.
Thus, we may apply the Poincar\'e recurrence theorem \cite{Walters1982} to right translation
by $g^{-1}$,
$$
R_{g^{-1}}\colon \Gamma\backslash G\longrightarrow \Gamma\backslash G,
$$
which is a measure preserving map. We obtain the following lemma.

\begin{fact}\label{fact:recurrence}
Let $\Gamma$ be a lattice in a Lie group $G$. Then, for every $g\in G$
and every nonempty open neighborhood $U$ of the identity, one has
$$
    \Gamma\cap Ug^nU^{-1}\neq\varnothing
$$
for infinitely many integers $n\geq 1$.
\end{fact}

\begin{proof}[Proof of Lemma \ref{lemma:dynamic}]
Fix $x\in\Omega$ and $h\in G_x$. Choose a relatively compact neighborhood
$U$ of the identity in $G$ such that
$$
    C:=\overline{U}\,x\subset\Omega.
$$
By Fact \ref{fact:recurrence}, for infinitely many $n$ there exist
$\gamma_n\in\Gamma$ and $u_n,v_n\in U$ such that
$\gamma_n=u_nh^nv_n^{-1}.$ Since $hx=x$, we have $\gamma_n(v_nx)=u_nx.$
Hence $ \gamma_n C\cap C\neq\varnothing.$
Since the action of $\Gamma$ on $\Omega$ is properly discontinuous and
$C$ is compact, only finitely many elements $\gamma\in\Gamma$ satisfy
$\gamma C\cap C\neq\varnothing$. Therefore, after passing to a subsequence,
we may assume that $\gamma_n=\gamma$ is constant. It follows that
$$
    h^n\in \overline{U}^{-1}\gamma\,\overline{U},
$$
which is compact. Hence $\{h^n\}_{n\in\ \N}$ admits a convergent subsequence.
\end{proof}

The next lemma combined with Lemma \ref{lemma:dynamic} will result the proof of  Proposition \ref{prop:dynamic}.

\begin{lemma}\label{lem: compacteness lemma}
Let $G$ be a connected Lie group. If every one-parameter
subgroup of $G$ has compact closure, then $G$ is compact.
\end{lemma}

The proof uses the following result of D'Ambra \cite{dambra}.

\begin{theorem}\cite[Theorem 7.0]{dambra}
Let $G$ be a connected Lie group. If every maximal
connected abelian subgroup of $G$ is compact, then
$G$ is compact.
\end{theorem}

\begin{proof}[Proof of Lemma \ref{lem: compacteness lemma}]
Let $A$ be a maximal connected abelian subgroup of $G$.
Since its closure is connected and abelian, $A$ is closed.
Thus, $A\simeq \mathbb T^k\times\mathbb R^\ell.$ If $\ell>0$, then $A$ contains a closed noncompact one-parameter subgroup of $G$, contradicting our assumption. Hence $\ell=0$, so $A$ is compact.
The conclusion follows from D'Ambra's theorem.
\end{proof}

\begin{lemma}
Let $H\leq G_x^\circ$ be a one-parameter subgroup.
Then $H$ has compact closure in $G_x^\circ$.
\end{lemma}

\begin{proof}
Let $h\in H$ be a nontrivial element. The closure of $H$
is either isomorphic to $\R$ or a torus $\mathbb T^k$
for some $k\geq 1$. Suppose, for a contradiction, that $\overline H\simeq\R$. Then the sequence $(h^n)_{n\geq 1}$ escapes every compact
subset of $G$. However, by Lemma \ref{lemma:dynamic}, this sequence
admits a convergent subsequence, a contradiction.
Hence $\overline H$ is compact.
\end{proof}

\subsection{Proof of Proposition \ref{prop:dynamic} and Theorem \ref{thm:solvable-complete}}

We have now all the tools to prove Proposition \ref{prop:dynamic} and so Theorem \ref{thm:solvable-complete}.

\begin{proof}[Proof of Proposition \ref{prop:dynamic}]
    Let $H$ be a one parameter group of $G_x^\circ$. Let $h\in H$ be a non-trivial element. Then we know by Proposition \ref{prop:dynamic} that $\{h^n\}_{n\in \mathbb{N}}$ is convergent, thus the closure of $H$ is compact. Thus by Corollary \ref{lem: compacteness lemma} we get that $G_x^\circ$ is compact.
\end{proof}

\begin{proof}[Proof of Theorem \ref{thm:solvable-complete}]
Since $\Gamma$ is finitely generated and linear, we may pass to a torsion-free subgroup
of finite index without affecting the conclusion. We claim that $\operatorname{cd}(\Gamma)\leq d$, where
$\operatorname{cd}(\Gamma)$ denotes the cohomological dimension
of $\Gamma$.

Since $\Gamma$ is solvable, by \cite{Crysta_FGH},
it admits a syndetic hull $S$ in
$\SL_d(\mathbb R)\ltimes\mathbb R^d$.
Thus, $S$ is a closed connected solvable Lie subgroup
containing $\Gamma$ as a uniform lattice.

Let $K_S$ be a maximal compact subgroup of $S$.
Fix $x\in\Omega$, and denote by $S_x$ its stabilizer
in $S$. By Proposition \ref{prop:dynamic},
the identity component $S_x^\circ$ is compact
and hence conjugate to a subgroup of $K_S$.
Therefore,
$$
\operatorname{cd}\Gamma
=\dim S-\dim K_S
\leq\dim S-\dim S_x^\circ
=\dim S-\dim S_x.
$$
However, $S/S_x$ injects in $\mathbb{R}^d$, hence $\dim S-\dim S_x\leq d.$ Then completeness follows from \cite[Theorem 3.5]{goldman1986affine}.
\end{proof}

\section{The rank-one case}\label{sec:rank1}

\subsection{General reduction result}\label{Section: geometric reuction}

The goal of this section is to prove Proposition \ref{mainprop}. For a semisimple Lie subgroup $G \leq \SL_d(\mathbb{R})$, we write $G = KA^+K$ for a chosen Cartan decomposition of $G$. For a subgroup $G' \leq G$, we denote by $A^+(G')$ the set of Cartan projections of elements of $G'$ in this decomposition.

Fix once and for all a Euclidean norm $|\cdot|$ on $\mathbb{R}^d$ that is invariant under the maximal compact subgroup $K$. 
We say that a subset $X$ of $G$ \emph{does not contract} $F$ with respect to $|\cdot|$ if $F$ is invariant under all elements of $X$ and if $|g \cdot x| \ge |x|$ for all $g \in X$ and all $x \in F$. 
In what follows, we will simply say that $X$ \emph{does not contract} $F$.

The proof of Proposition \ref{mainprop} requires some preparation and is inspired by arguments from \cite{tholozan2015completude}. To fix notation, let $p \in \llbracket 1, d \rrbracket$, and consider $\operatorname{Gr}_p(\mathbb{R}^d)$, the Grassmannian of $p$-dimensional subspaces, note that the linear action of $\mathrm{SL}_d(\mathbb{R})$ induces a natural action on $\operatorname{Gr}_p(\mathbb{R}^d)$. We say that $E, F\in \operatorname{Gr}_p(\mathbb{R}^d)$ are \emph{transverse} if $E + F = \mathbb{R}^d$. Throughout this section, we work under the assumptions of Proposition \ref{mainprop}. Write $\Omega = \Gamma \cdot \mathcal{C}$ with $\mathcal{C} \subset \Omega$ compact. Let $\varepsilon > 0$ be such that $\mathcal{C} + B_\varepsilon \subset \Omega$, where $B_\varepsilon$ denotes the Euclidean ball of radius $\varepsilon$ centered at $0$, defined with respect to the Euclidean norm $\lvert \cdot \rvert$.

\begin{lemma}\label{prop_variétés_discompacité_supérieure}
Let $y \in \partial \Omega$. Then there exist sequences $(g_n)_{n\in\mathbb{N}} \subset \Gamma$, $(x_n)_{n\in\mathbb{N}} \subset \mathcal{C}$, and $(E_n)_{n\in\mathbb{N}} \subset \mathcal{F}$ such that, up to extraction of a subsequence:
\begin{enumerate}
    \item $y_n:=g_n \cdot x_n$ converges to $y$ as $n \to \infty$;
    \item for every $n$, one has $(E_n+y_n) \cap B(y_n,\varepsilon) \subset \Omega$;
    \item $E_n$ converges in $\operatorname{Gr}_p(\mathbb{R}^d)$ to some element $E_y$ in the closure of $\mathcal{F}$, with $(E_y+y) \cap B(y,\varepsilon)\subset \partial \Omega$, where $B(y,\varepsilon)$ is Euclidean ball of radius $\varepsilon$ centered at $y$.
\end{enumerate}
\end{lemma}

\begin{proof}
Since $\Omega = \Gamma \cdot \mathcal{C}$, there exist sequences $(x_n)_{n\in \mathbb{N}} \subset \mathcal{C}$ and $(g_n)_{n\in \mathbb{N}} \subset \Gamma$ such that $g_n \cdot x_n \to y$. This proves the first point. For all $n \in \mathbb{N}$, we have
\begin{equation}\label{eq:lemma3.3-1}
    g_n \cdot (x_n + B_\varepsilon) = y_n + L(g_n) \cdot B_\varepsilon \subset \Omega.
\end{equation}
For each $n$, write the Cartan decomposition $L(g_n) = w_n a_n w_n'$, with $w_n, w_n' \in K$ and $a_n \in A^+(L(\Gamma))$. Let $E$ be the $p$-subspace appearing in the hypothesis of Proposition \ref{mainprop}, and set $E_n := w_n \cdot E$. We claim that $L(g_n)\cdot B_{\varepsilon}$ contains $E_n \cap B_{\varepsilon}$. Indeed, let $v \in w_n \cdot E \cap B_\varepsilon$ and set $x = w_n^{-1}v \in E$. Since $K$ preserves the Euclidean norm $\lvert\cdot\rvert$ and $a_n$ does not contract $E$, we have $\lvert a_n^{-1}x\rvert \le \lvert x\rvert = \lvert v\rvert \le \varepsilon$. Therefore $u := w_n'^{-1}(a_n^{-1}x)$ belongs to $B_\varepsilon$, and $L(g_n)u = w_n a_n w_n' \bigl(w_n'^{-1}a_n^{-1}x\bigr) = w_n x = v$. Hence $w_n \cdot E \cap B_\varepsilon \subset L(g_n)\cdot B_\varepsilon$, which proves the claim.

\smallskip

\noindent
It follows from \eqref{eq:lemma3.3-1} that $\Omega$ contains $y_n + (E_n \cap B_\varepsilon)$, establishing the second point. By compactness of $\overline{\mathcal F}$, we may assume (up to extracting a subsequence) that $E_n$ converges in the Grassmannian to some $E_y \in \overline{\mathcal F}$. Then $\overline{\Omega}$ contains $y + (E_y \cap B_\varepsilon)$, and since $y \in \partial \Omega$ and the $\Gamma$ action on $\Omega$ is proper, the inclusion must in fact be
$$
    (E_y \cap B_\varepsilon) + y \subset \partial \Omega.
$$
This proves the third point and completes the proof.
\end{proof}
\begin{remark}
    In the previous lemma, we do not use the transversality of the elements of the closure of $\mathcal{F}$. In particular, the lemma asserts that each point of the boundary contains a small neighborhood of an affine subspace.
\end{remark}

The next step is to show the uniqueness of the subspace $E_y$. Before doing so, we record the following elementary fact about transversality.

\begin{lemma}\cite[Lemma 3.8]{diaf2026completeness}\label{transversal_intersection}
Let $E,F\subset\mathbb{R}^d$ be linear subspaces with $E+F=\mathbb{R}^d$.
Let $E_n\to E$ in the Grassmannian and $y_n\to y$ in $\mathbb{R}^d$.
Fix $\varepsilon>0$. Then for all sufficiently large $n$ the sets
$$
(E_n+y_n)\cap B(y_n,\varepsilon)\qquad\text{and}\qquad (F+y)\cap B(y,\varepsilon)
$$
have nonempty intersection.
\end{lemma}
We arrive at the following lemma.
\begin{lemma}\label{uniqueness_of_plane}
Let $y \in \partial \Omega$, and let $E_y$ be an element of the closure of $\mathcal{F}$ as in Lemma \ref{prop_variétés_discompacité_supérieure}. Then $E_y$ is the unique element of $\overline{\mathcal{F}}$ with the property that there exists a neighborhood $\mathcal{V}$ of $y$ in $\mathbb{R}^d$ such that $(y + E_y) \cap \mathcal{V} \subset \partial \Omega$.
\end{lemma}

\begin{proof}
Let $(y_n)_{n \in \mathbb{N}}$, $(E_n)_{n \in \mathbb{N}}$, and $E_y$ be as in Lemma \ref{prop_variétés_discompacité_supérieure}. By shrinking $\varepsilon$ if necessary, we may assume that $B(y, \varepsilon) \subset \mathcal{V}$.

We proceed by contradiction and assume the existence of a subspace $E \in \overline{\mathcal{F}}$, distinct from $E_y$, that satisfies the hypothesis of the lemma. Then $E_y$ and $E$ are transverse. Since $y_n \to y$ and $E_n \to E_y$, and using Lemma \ref{transversal_intersection}, we deduce that, for $n$ large enough, the set $(E_n + y_n) \cap B(y_n, \varepsilon)$ intersects $(E + y) \cap B(y, \varepsilon)$.

However, by Lemma \ref{prop_variétés_discompacité_supérieure}, we have $(E_n + y_n) \cap B(y_n, \varepsilon) \subset \Omega$, while $(E + y) \cap B(y, \varepsilon) \subset (E + y) \cap \mathcal{V} \subset \partial \Omega$, a contradiction.
\end{proof}

The next lemma goes a step further by showing that $\partial \Omega$ contains an affine subspace.

\begin{lemma}\label{lem_H_inclus_partial_O}
Let $y \in \partial \Omega$, and let $E_y$ be as in Proposition \ref{prop_variétés_discompacité_supérieure}. Then $E_y + y \subset \partial \Omega$.
\end{lemma}

\begin{proof}
Let $\mathcal{U}$ be the set of points $z \in E_y + y$ such that $(E_y + y) \cap B(z, \varepsilon) \subset \partial \Omega$. We will show that $\mathcal{U}$ is both open and closed (and nonempty) in $E_y + y$, so that $\mathcal{U} = E_y + y$, which will complete the proof.

Clearly, the set $\mathcal{U}$ is closed and contains $y$. We now show that $\mathcal{U}$ is also open. Let $z \in \mathcal{U}$ and $z' \in (E_y + y) \cap B(z, \varepsilon / 2)$. Then $(E_y + y) \cap B(z', \varepsilon / 2) \subset \partial \Omega$. Thus, by Lemma \ref{uniqueness_of_plane}, we have $E_y = E_{z'}$.
By the definition of $E_{z'}$, this implies
$$(z'+E_y)\cap B(z',\varepsilon)=(z'+E_{z'})\cap B(z',\varepsilon)
\subset \partial\Omega.$$
Since $z'\in y+E_y$, we have $z'+E_y=y+E_y$.
Hence $z'\in\mathcal U$. Therefore $\mathcal U$ is open,
which completes the proof.
\end{proof}
We now have all the tools to prove Proposition \ref{mainprop}.

\begin{proof}[Proof of Proposition \ref{mainprop}]
First, we claim that two affine subspaces contained in $\partial \Omega$ whose linear parts lie in $\overline{\mathcal{F}}$ must be parallel. By Lemma \ref{lem_H_inclus_partial_O}, there exists at least one affine subspace contained in $\partial \Omega$ whose linear part belongs to $\overline{\mathcal{F}}$.

Assume, for contradiction, that there exist two affine subspaces $E$ and $E'$ contained in $\partial \Omega$ with distinct linear parts in $\overline{\mathcal{F}}$. Then $E$ and $E'$ are transverse, so they intersect at some point $y \in \partial \Omega$, which contradicts Lemma \ref{uniqueness_of_plane}. Therefore $E$ and $E'$ must have the same linear part; in other words $E$ and $E'$ are parallel.

Next, let $E_0$ denote the common linear part of the hyperplanes contained in $\partial \Omega$. Since $\Gamma$ preserves $\Omega$, the group $L(\Gamma)$ preserves $E_0$. Let $x \in \Omega$ and suppose that $(E_0 + x)$ is not contained in $\Omega$. Then there exists $y \in (x + E_0) \cap \partial \Omega$, and by the previous argument $(y + E_0) \subset \partial \Omega$. Hence $x \in \partial \Omega$, a contradiction. This completes the proof.
\end{proof}

\subsection{Transversality phenomenon for rank 1 Lie groups}

Let $G$ be a connected real semisimple Lie group of real rank one, and let
$$
\rho:G\longrightarrow \mathrm{GL}(V)
$$
be a faithful finite-dimensional real representation.
The goal of this section is to apply the general reduction result of Proposition \ref{mainprop} to a case of interest. We show  the following:
\begin{prop}\label{thm_rank_one}
Let $G$ be a connected real semisimple Lie group of real
rank one, and let $\rho:G\longrightarrow\SL_d(\mathbb R)$ be a faithful representation. Let $\Omega\subsetneq\mathbb R^d$ be a domain, and let
$\Gamma\leq\rho(G)\ltimes\mathbb R^d$ be a discrete subgroup
acting freely, properly discontinuously, and cocompactly
on $\Omega$.

Then there exists a linear subspace
$E\subset\mathbb R^d$, with $\dim E\geq d/2$,
which is invariant under $L(\Gamma)$, such that
$\Omega$ is foliated by affine subspaces parallel to $E$.
Moreover, the stabilizer of $E$ in $G$ is a   
parabolic subgroup $P'$ of $G$. In particular, $L(\Gamma)\leq \rho(P').$
\end{prop}

For that we need some algebraic preparation in order to apply Proposition \ref{mainprop}. Fix a Cartan decomposition $\g=\k\oplus\p$, and let $K$ be the
corresponding maximal compact subgroup of $G$.
Let $\a\subset\p$ be a maximal abelian subspace. Since $G$ has
real rank one, $\dim\a=1$. Choose a nonzero element $H\in\a$, and set
$$
\a^+=\mathbb{R}_{\geq0}H,
\qquad
A=\exp(\a),
\qquad
A^+=\exp(\a^+).
$$
Let $\Sigma=\Sigma(\g,\a)$ be the restricted-root system, and
denote by
$$
\Sigma^+=\{\alpha\in\Sigma:\alpha(H)>0\}
$$
the corresponding set of positive restricted roots.
For $\alpha\in\Sigma$, let
$$
\g_\alpha
=\{X\in\g:[Y,X]=\alpha(Y)X
\text{ for all }Y\in\a\}.
$$
Set
$$
\n=\bigoplus_{\alpha\in\Sigma^+}\g_\alpha,
\qquad
N=\exp(\n),
\qquad
M=Z_K(A),
\qquad
P=MAN.
$$
so that $P$ is a parabolic subgroup of $G$. Recall that since~$G$ has rank one,~$P$ is the only parabolic subgroup of~$G$ up to conjugation.

Since $G$ has real rank one, its restricted Weyl group
$W=N_K(A)/M$ has order two.
Let $w\in N_K(A)$ represent its nontrivial element, so that
\begin{equation}\label{eq:weyl}
\operatorname{Ad}(w)H=-H.
\end{equation}
We shall use the following standard decompositions:
\begin{align}
G&=KAN=KP,
\qquad K\cap P=M,
&&\text{(Iwasawa)}\label{eq:iwasawa}\\
G&=P\sqcup NwP,
&&\text{(Bruhat)}\label{eq:bruhat}\\
G&=KA^+K.
&&\text{(Cartan)}\label{eq:cartan}
\end{align}
Let
$$
d\rho:\g\longrightarrow\mathfrak{gl}(V)
$$
be the differential of $\rho$ at the identity.
By a standard result for finite-dimensional representations
of semisimple Lie groups, there exists an inner product on $V$
for which $d\rho(X)$ is skew-symmetric for every $X\in\k$
and symmetric for every $X\in\p$.
We denote by $\|\cdot\|$ the associated Euclidean norm,
which is invariant under $\rho(K)$.
In particular, $d\rho(H)$ is symmetric and hence
diagonalizable over $\mathbb R$.

Let $\Sigma(V,\a)\subset\a^*$ denote the set of weights of
the representation $d\rho|_{\a}$. We have the decomposition
$$
V=\bigoplus_{\lambda\in\Sigma(V,\a)}V_\lambda,
$$
where
$$
V_\lambda
=\{v\in V:d\rho(Y)v=\lambda(Y)v
\text{ for all }Y\in\a\}.
$$
Define
$$
V_+
=\bigoplus_{\lambda(H)\geq0}V_\lambda,
\qquad
V_-
=\bigoplus_{\lambda(H)\leq0}V_\lambda,
\qquad
V_0=V_+\cap V_-.
$$

\begin{lemma}\label{lem:P}
The stabilizer of $V_+$ in $G$ is the   
parabolic subgroup $P$:
$$
\operatorname{Stab}_G(V_+)=P.
$$
\end{lemma}

\begin{proof}
Set $Q=\operatorname{Stab}_G(V_+)$. We first show that $P\subseteq Q$.
For every $\alpha\in\Sigma^+$ and $X\in\g_\alpha$, it is not difficult to check the following weight-shift identity
$$
d\rho(X)V_\lambda\subseteq V_{\lambda+\alpha}.
$$
Since $\alpha(H)>0$ for every $\alpha\in \Sigma^+$, we deduce that the subgroup $N$ preserves
$V_+$. Moreover, both $A$ and $M$ preserve each
weight space $V_\lambda$. Thus $P=MAN$ preserves
$V_+$, proving that $P\subseteq Q$.

Suppose now that $Q\neq P$.
By the Bruhat decomposition \eqref{eq:bruhat}, there exists an element
$q=nwp\in Q$, with $n\in N$ and $p\in P$.
Since $N\subseteq P\subseteq Q$, it follows that $w\in Q$.

On the other hand, $\operatorname{Ad}(w)H=-H$ implies $\rho(w)V_\lambda=V_{-\lambda},$ and consequently
\begin{equation}\label{eq:reverse}
\rho(w)V_+=V_-.
\end{equation}
Since $w\in Q$, we obtain $V_+=V_-.$ This forces every weight of $d\rho|_{\a}$ to vanish,
so $d\rho(H)=0$. As $H\neq0$, this contradicts
the faithfulness of $\rho$. Therefore $Q=P$, as claimed.
\end{proof}

Let $k=\dim V_+$. Since $P$ stabilizes $V_+$ by
Lemma \ref{lem:P}, the following map is well-defined:
$$
\Phi:G/P\longrightarrow\operatorname{Gr}_k(V),
\qquad
\Phi(gP)=\rho(g)V_+.
$$

\begin{lemma}\label{lemma:transversality rank 1}
The map $\Phi$ is $G$-equivariant, and injective, and its image is
\begin{equation}\label{eq:my orbit}
 \mathcal F=\{\rho(k)V_+:k\in K\}\cong K/M.    
\end{equation}
For distinct $x,y\in G/P$, we have that
$$
 \Phi(x)+\Phi(y)=V
$$
\end{lemma}

\begin{proof} Equivariance is immediate.
By Lemma \ref{lem:P}, we have $\Stab_G(V_+)=P$, so $\Phi$ is
injective. Moreover, the Iwasawa decomposition \eqref{eq:iwasawa}
gives $G=KP$ and $K\cap P=M$. Hence,
$$
\mathcal F=\{\rho(k)V_+:k\in K\}\cong K/M.
$$

For distinct $x=g_1P$ and $y=g_2P$, write $g_1^{-1}g_2=nwp$ by Bruhat decomposition \eqref{eq:bruhat}, with $n\in N$ and $p\in P$. Lemmas \ref{lem:P}
and the fact $\rho(w)V_+=V_-$ (see \eqref{eq:reverse}) imply
$$
 \Phi(x)=\rho(g_1n)V_+,\qquad
 \Phi(y)=\rho(g_1n)V_-.
$$
Since
$V_++V_-=V$, applying the invertible map $\rho(g_1n)$
shows the claim.
\end{proof}

\begin{proof}[Proof of Proposition \ref{thm_rank_one}]
Since $G$ is connected and semisimple, the image of
$\rho$ is contained in $\SL(V)$. Moreover, $\rho(G)$
is a closed semisimple Lie subgroup of $\SL(V)$.
We show that the family $\mathcal F$ constructed in
Lemma \ref{lemma:transversality rank 1} satisfies the
hypotheses of Proposition \ref{mainprop}.

By the Iwasawa decomposition $G=KP$ and the fact that
$P$ stabilizes $V_+$, we deduce that $\mathcal F$ is a
$\rho(G)$-orbit in $\operatorname{Gr}_k(V)$.
Moreover, $\mathcal F$ is compact, since it is also
a $\rho(K)$-orbit and $K$ is compact. Moreover, by
Lemma \ref{lemma:transversality rank 1},
any two distinct elements of $\mathcal F$ are transverse.

It remains to verify the non-contracting condition.
Recall that $V$ is equipped with a $\rho(K)$-invariant
Euclidean norm $\|\cdot\|$ for which $d\rho(H)$ is symmetric.
In particular, the weight spaces are mutually orthogonal.
For every $v=\sum_{\lambda(H)\geq0}v_\lambda\in V_+$
and $t\geq0$, we have
$$
\|\rho(\exp(tH))v\|^2=\sum_{\lambda(H)\geq0}
e^{2t\lambda(H)}\|v_\lambda\|^2
\geq\sum_{\lambda(H)\geq0}\|v_\lambda\|^2
=\|v\|^2.
$$
Thus $\rho(A^+)$ does not contract $V_+$.
Since the Cartan decomposition of $\rho(G)$ is
$$
\rho(G)=\rho(K)\rho(A^+)\rho(K),
$$
the Cartan projections of $L(\Gamma)$ do not
contract $V_+$ either.

Therefore, Proposition \ref{mainprop} applies and
provides an $L(\Gamma)$-invariant subspace
$E\in\mathcal F$ such that $\Omega$ is foliated
by affine subspaces parallel to $E$ and by Lemma \ref{lem:P} $L(\Gamma)$ is contained is a conjugate of the    parabolic subgroup $P$. This completes the proof.
\end{proof}

\subsection{Proof of Theorem \ref{thm:rank1-complete}}
We now prove Theorem \ref{thm:rank1-complete} using the
geometric reduction established in Proposition \ref{thm_rank_one}
and the amenability of    parabolic subgroups.
\begin{proof}[Proof of Theorem \ref{thm:rank1-complete}]
Suppose, by contradiction, that $\Omega\neq\mathbb R^d$.
By Proposition \ref{thm_rank_one}, there exists a   
parabolic subgroup $P$ of $G$ such that $\Gamma\leq\rho(P)\ltimes\mathbb R^d.$ Since $P$ is amenable, so is
$\rho(P)\ltimes\mathbb R^d$, which is a Lie group
with finitely many connected components.
By a result of Milnor \cite[Lemma 2.2]{milnor1977fundamental},
the discrete subgroup $\Gamma$ is virtually polycyclic,
and hence virtually solvable. Then the proof follows by Theorem \ref{thm:solvable-complete}.
\end{proof}

\bibliographystyle{alpha}
\bibliography{bib}

\end{document}